\documentclass[8pt,twocolumn]{article}
\usepackage[dvips]{graphicx}
\usepackage{amsmath}
\usepackage{amssymb}
\usepackage{array}

\begin{document}

\bf\normalsize\title{Symmetry in Partial Sums of $n^{-s}$ \\\
\textit{(\textnormal{\textbf{or}}, Critical Line Zeros of $\zeta$ Are Easy)}}
\author{George H. Nickel}

\maketitle
\normalsize

\section{Key Words}

Riemann $\boldsymbol{\zeta}$-function, analytic continuation, functional equation, $\boldsymbol{\zeta}$-zeros.

\rm

\section{Summary}

The Riemann $\zeta$-function is $\sum_1^{\infty}n^{-s}$ for $\sigma >1$, or its analytic continuation for $\sigma < 1$. It was originally studied by Euler on the integers and extended to the complex plane by Riemann for the purpose of investigating the properties of prime numbers.  Today, it finds application in statistical physics and quantum field theory\cite{wolf}.

An inherent symmetry underlies the dependence of the partial sums $\sum_1^n m^{-s}$ on $n$. This symmetry relates individual \textit{steps}\footnote{$[x]$ denotes the integer floor of $x$} $n<n_p = [\sqrt{t/2\pi}]$ to collections of adjacent steps $n > n_p$.  It defines two particular points in the complex plane: \\

\begin{quote} A \textit{center point} conjugate to itself in a sense described below, with an axis of bilateral symmetry passing through it at a specific angle.\end{quote}

 \begin{quote}A point conjugate to the origin, beyond which the sum of steps either converges or approaches infinity, depending on the value of $\sigma$.\end{quote}

 Forward step angle differences $\delta\theta_n = -t/(n+1/2) + \ldots$, $\delta^2\theta_n = t/(n+1/2)^2 + \ldots$, etc. are defined between adjacent steps.  At the point referred to here as $n_p$, defined by the partial sum $\sum_1^{n_p}m^{-s}$, the second angle difference $\delta^2 \theta_{n_p}$ equals $2\pi$, equivalent to zero. As $n$ increase beyond $n_p$, the second angle variation $\delta^2\theta_n$ decreases monotonically from $2\pi$ to $0$, while first angle variation $\delta\theta_n$ decreases through $n_p$ regions bounded by steps whose first angle differences are odd multiples of $\pi$.  Each of the $n_p$ regions $[n_p^2/(n+1)]<N<[n_p^2/2\pi (n-1)]$ is conjugate to the single initial step $n$.  There are approximately $n_p^2/n^2$ steps in the region conjugate to $n$. Because of truncation to step numbers $n \in \mathbf{P}$, the precise boundaries of conjugate regions are not exactly coincident with centers of steps.

 When $n_p >> 1$, angles and lengths of the steps immediately adjacent to $n_p$ are nearly constant, forming a polygon tangent to step $n_p$. $P(s)$, the center of symmetry, is recognized here to be the center of this feature, laterally displaced from $n_p$.  The specific form of the conjugacy is that each initial step $n^{-s}$, transformed to $n^{s-1}$ and multiplied by a complex factor independent of $n$, equals the sum of steps in its conjugate region. The polar form of the complex factor is $Q(s)=(t/2\pi)^{1/2-\sigma}e^{2i\Theta}$, where $\Theta(s)$ is approximately the sum of the angles of steps $n_p$ and $n_{p+1}$.  To the next order of approximation, end corrections $\Delta_n$ for the conjugate regions give the precise points conjugate to the ends of step $n$, assumed to correspond to the geometric center of the constellation of steps having $\delta \theta \simeq$ the $nth$ odd multiple of $\pi$ at $N_n$.  The lengths of a conjugate region and its initial step $n$ are equal when $\vert Q \vert = 1$, which only occurs for $\sigma = 1/2$. $Q(s) \ne 0$ on the critical strip. The sum of the angles of each initial step and its summed conjugate region is $2\Theta(s)$ and the axis of symmetry is at the angle $\Theta \pm \pi/2$.  Angle differences of adjacent regions conjugate to $n,n+1$ can be summed immediately, and are seen to equal $\delta\theta_n$. The sum of steps in a conjugate region requires explicit summation, which can be approximated by a Jacobi $\theta$-function for a large number of conjugate steps.  This is the case for $t \gg 1$ and $n \ll n_p$.

 The description of $\zeta(s)$ in terms of symmetry is shown in subsequent sections to agree with analytic results.  The continuation of $\zeta(s)$ for $\sigma < 1$ corresponds to the point $O'(s)$ conjugate to the origin, using Riemann's general definition. As a specific example, $O' \to \zeta(s)$ for $\sigma < 1$ is consistent with the Euler-Maclaurin procedure.  The functional equation $\zeta(s) = Q(s) \zeta(1-s)$ follows from the symmetry.  The central symmetry point $P(s)$, which is conjugate to itself, is related directly the Riemann-Siegel equation; the lowest order geometric construction reproduces Riemann's first order of correction. Other specific results obtained by Riemann and his successors are addressed, including the number of zeros of $\zeta(\sigma+it)$ for $t \le T$ as a function of $\Theta(s)$ and their relation to the critical line $\sigma=1/2$.

In this work, detailed derivations and figures are omitted in favor of brevity. A more complete discussion, including bounds, generalized $\zeta$ functions and the special significance of prime steps, is found in \cite{bulletin}.

\rm
\section{A Review of the Traditional Analysis}

Mathematicians of the early 19th century obtained general expressions for functions previously defined only on the integers, and extended them to the complex plane.\footnote{See, for example, \cite[Chap.~12]{apostol}. These topics are chosen both to illustrate the history and to introduce specific functions encountered below.} Euler derived an expression for the factorial in terms of the definite integral familiar to all mathematics students

\begin{equation}
\Gamma(x) = \int_0^{\infty} y^x e^{-y} dy \to (x-1)! \quad for \quad x \in \mathbb{P}
\end{equation}

\noindent which Gauss expressed alternately as

\begin{equation}
\Gamma(x) = \lim_{N \rightarrow \infty} \frac{1 \cdot 2 \cdot \ldots (N-1)}{(x+1)(x+2) \ldots (x+N-1)}N^x
\end{equation}

Euler's integral arises naturally in analysis, and Gauss' form is useful to define connections between $\Gamma$-functions. For example, the functional equation $\Gamma (x+1) = x \Gamma(x)$ is supplemented by $sin\pi s =\pi s/\Gamma(s)\Gamma(1-s)$ and $\Gamma(s)=2\Gamma((s+1)/2)\Gamma(s/2)$ .

A discrete series can often be represented by the expansion of a continuous \textit{generating function}\cite{wilf}; this facilitates integrations in the complex plane. An elementary application from statistical physics is the integral of the Planck distribution for photons\footnote{The numerator $x^3 dx$ arises from the dimensionless form of the photon energy $h \nu$ times $(h\nu /c)^2 d(h\nu /c)$ for a volume element of momentum space, and the denominator represents the quantum distribution function for bosons, where $x=h\nu/kT$.}

\begin{equation}
\int_0^{\infty} \frac{x^3 dx}{e^x-1} dx,
\end{equation}

Expansion of the integrand gives the series under discussion here with argument $s=4$:

\begin{equation}
\int_0^{\infty} x^3 \sum_1^{\infty} e^{nx} dx = \Gamma (4) \sum_1^{\infty} n^{-4} \to \Gamma (4) \zeta(4)
\end{equation}

 Riemann generalized this to the complex plane by replacing the numerator $x^3$ of this example by $(-x)^s$ and using a contour from $\infty$ above the positive real axis, circling the origin in the positive sense at radius $\delta$, followed by a return to $\infty$ below the positive real axis. He showed that this procedure yields analytic continuation for $\sigma < 1$\cite[Sec.~1.4]{edwards}.

Euler had earlier expressed the fundamental Dirichlet series as a product of terms involving the prime numbers (his ``Golden Rule'')

\begin{equation}
\sum_1^{\infty} \frac {1}{n^{x}} = \prod_{p} \frac{1}{(1-\frac{1}{p^x})}.
\end{equation}

Riemann used the poles resulting from the logarithm of this product as the starting point of his analysis leading to a relation between zeros of $\zeta$ and the primes, as well as an expression for $\Pi(x)$, the number of primes less than or equal to $x$ \cite[Sec.~1.18]{edwards}

\begin{equation}
\Pi(x) = Li(x)-\sum_{zeros \rho} Li(x^{\rho})+\ldots
\end{equation}

His first proof of the \textit{functional equation of the $\zeta$-function} under the $s \to 1-s$ transformation \cite[Sec.~1.6]{edwards}

\begin{equation}
\zeta (s) = \frac{(-2\pi i)^s}{\Gamma (s)} \zeta(1-s)
\label{funct}
\end{equation}

\noindent used the contour integration above, in addition to disjoint circles about the poles of the generating function on the imaginary axis in the negative sense.  The two contributions add to zero, giving Equation \ref{funct}.  It will be seen that this is a detailed relation, applying not only globally to the sum over all steps and poles, but to the individual pairs of terms as well.

Jacobi's $\theta$-functions include

\begin{equation}
G(u)=\sum_{-\infty}^{\infty} e^{-n^2 \pi u^2} \quad and \quad \Psi(u)=\sum_1^{\infty} e^{-n^2 \pi u}
\end{equation}

\noindent which Riemann used for the second proof of his functional equation \cite[Sec.~1.7]{edwards}. He then applied Gauss' $\Gamma$-function identities to express a function $\xi(s)$ which is invariant under $s \to 1-s$ as

\begin{equation}
\xi(x)=(s-1)\Gamma(1+s/2) \pi^{-s/2}\zeta(s)
\end{equation}

The history of the use of these integrations and inversion transforms to examine the prime numbers is described in detail by Edwards \cite{edwards}. General methods of analytic number theory are discussed by Apostol \cite{apostol}.

\section{Angles and Conjugate Regions}

\begin{figure*}[pt]
\centering
\includegraphics[angle=0,width=.7\textwidth]{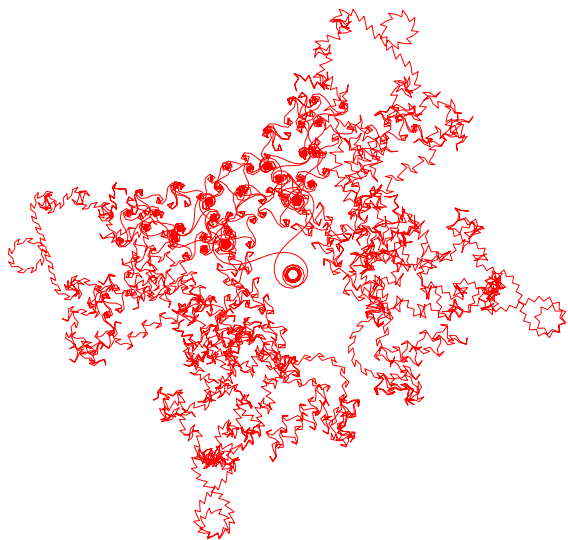}
\caption{Steps near the symmetry center for $s = 1/2+10^9 i$}
\label{blob1pt9}
\end{figure*}

For extremely large argument $t$ and $n \ll t/\pi$, the first angle differences contain many multiples of $2\pi$.  Although step angles vary monotonically as $-tlog(n)$, their modulus with respect to $2\pi$ appears to vary irregularly when $\delta \theta \gg 2\pi$.  For one particular number $n_p$ of initial steps, there are $n_p$ subsequent steps having $\delta \theta$ equal to an even multiple of $\pi$.  Solving $n_p = \delta \theta_{n_p}/2\pi$ gives $n_p =[\sqrt{t/2\pi}-1/2]$.\footnote{Although derived from a different approach, this quantity plays a significant role in Riemann's analysis. It can also be defined by the condition $\delta^2 \theta = 2\pi$.}  For $n_p \gg 1$ its adjacent steps approximate the ``pendant'', a regular polygon with slowly varying lengths and $\delta \theta$, a distinct feature in the step plot.  See Figure \ref{blob1pt9}. For all steps $n > n_p$, $\delta\theta_n$ decreases smoothly through odd and even multiples of $\pi$, giving $n_p$ regions conjugate to the $n_p$ initial steps. The center step of the region conjugate to step $n$ is an inflection point whose the first angle difference is an even multiple of $\pi$ and its end points are at the steps whose first angle differences are odd multiples of $\pi$.  This gives each conjugate region the characteristic ``f-slot'' shape on the belly of a stringed instrument, to be described in detail below.

\begin{figure*}[p]
\center
%\put(0,0)
\includegraphics[angle=0,width=.7\textwidth]{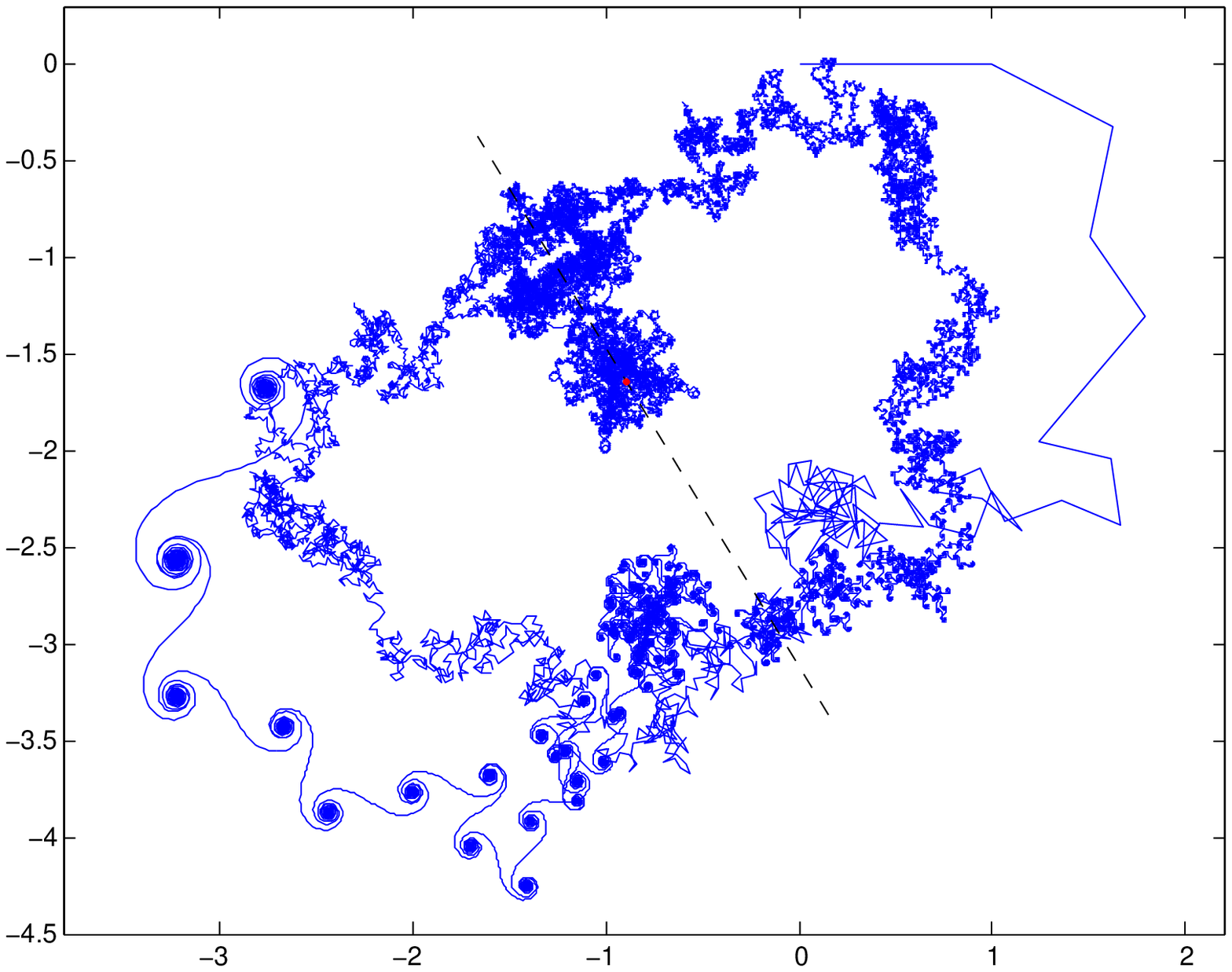}
\caption{Plot of summed steps for $\zeta(1/2+ 10^9i)$. }
\label{xy9}
\end{figure*}

 Figure \ref{xy9}, the plot of $t/\pi \simeq 300$ million steps for $s=1/2+10^9i$, illustrates these details.  For this large value of $t$, the angle difference regularities of alternating odd and even multiples of $\pi$ are clearly evident for steps and their conjugate regions. A red dot marks the center of symmetry of this step plot, which shows that the direction angles of each initial step and its conjugate region add to a constant, whose perpendicular defines an axis of symmetry (dashed in the figure). In particular, since the first step is at the angle $0$ with length $1$, the angle of the final conjugate region \textit{is} this constant angle sum, denoted here by $2 \Theta(s)$.

 A central proposition of this work is that the geometric center of each scroll corresponds precisely to the conjugate point of an initial step. Because of step discretization, this differs from the midpoint of a step by a longitudinal displacement $\lambda$ which depends on $\sigma$ and a transverse displacement $i\tau(t)$ from the step center.  A first-order derivation of these displacements is given below.

The density of steps near the symmetry center of this figure precludes visualization of its details; these are more clearly displayed in Figure \ref{blob1pt9}, which reveals other interesting features such as the central pendant with its self-symmetric center and fractal-like repetitions where angle changes reproduce smaller representations of the pendant. There are many conjugate figures having rational ratios of steps in their respective counterparts.  These details are described and illustrated in \cite{bulletin}.

\section{Extent of Conjugate Regions}

The extent of the region conjugate to step $n<n_p$ is given by

\begin{equation}
\sum_{N_n}^{N_{n-1}} m^{-s}, \quad N_{n}=[\frac{n_p^2}{n-1/2} -1/2].
\end{equation}

The inflection point of the region conjugate to step $n$ is at step $N_n^0 = [n_p^2/(n-1/2)]$, where the step angle is $\theta_{N_n^0} =-t log(N_n^0)$, the first angle difference $mod,2\pi$ is 0, and the second angle difference is $2\pi (n_p/N_n^0)^2$. An approximate result for the sum is given by Jacobi's $\theta$-function $G(u)$. This is appropriate for large $t$ and $n$ sufficiently less than $n_p$ that there are many steps in the sum.  The approximately linear variation of step length across the conjugate region can be neglected, since it sums to zero, and the second order variation is $\mathcal{O} (N_n^0)^{-2}$.  Similarly, the third angle difference is $\mathcal{O} (N_n^0)^{-3}$.  To this order of approximation, step angles can be assumed to vary as the square of $\delta n = n-N_n^0$.  This leads to the classic $\textit{Cornu}$ or $\textit{Euler}$ spiral, where the angle of the tangent to the curve varies quadratically with arc length.  The limits of the sum can be taken as $\pm \infty$ because the spirals near the end scrolls oscillate around a central value, giving

\begin{equation}
\sum_{N_n}^{N_{(n-1)}} m^{-s} = {N_n^0}^{-s}\sum_{-\infty}^{\infty} e^{ i \pi \delta n^2\frac{n_p^2}{{N_n^0}^2}}
\end{equation}

The approximation\footnote{Valid to five places for $ \lvert u \rvert = \lvert i^{1/2}\frac{N}{n_p} \rvert < 1/2$} $G(u) = \sum_{-\infty}^{\infty} e^{-\pi n^2 u^2} \simeq 1/u$ gives

\begin{equation}
\sum_{N_n}^{N_{n-1}} m^{-s} = \frac{n_p^{1-2 \sigma}}{n^{1-s}} e^{i 2\theta_{N_n^0} -i\pi/4 + \ldots}
\label{sumconj}
\end{equation}

\noindent where $\theta_{N_n^0}$ is the angle of the step at the symmetry center.  Algebraically, the angle $\pi /4$ is the term $\sqrt{i}$ in the definition of $u$; from Figure \ref{xy9} it is seen to be the angle between the central step about which the integral was evaluated and the summed extent of the Euler spiral.  A further angle term, denoted by $\delta t$, is due to the difference between $t$ and the exact value for which $\delta^2 \theta_{N_o} \equiv \pi N^2/n_p^2$.  Because the reference value of $\delta t$ is an integer multiple of $2\pi$, this correction is denoted in \cite{edwards} simply as $t$.  Higher order angle terms due to the approximations of this approach are bounded by the first omitted term from Riemann's analysis, $1/48t$.

The ratio of $\sum_{N_-}^{N_+} n^{-s}$  to $n^{1-s}$ is central to this work, as it quantifies the overall symmetry of the step plot.  It is denoted here by

\begin{equation}
Q(s) = n_p^{1-2 \sigma} e^{ 2i\theta_{n_p} -\pi/4+ t +\ldots}=n_p^{1-2 \sigma} e^{2i\Theta}
\end{equation}

\noindent which, except for the choice of sign, is identical to Riemann's result, as seen in the functional equation\cite[Sec.~6.5]{edwards}. It can be shown that $\zeta(s)$ lies on a line through the origin at the angle $\Theta =  \theta_{n_p} -\pi/8+ t/2 +\ldots$.

The factor $Q(s)$ can be related to the results of analytic number theory.  Consider the comparison of the integral around the $nth$ imaginary axis pole pair of the generating function to the contour integral of the single step $n$:

\begin{equation}
(2\pi i)(-2\pi i n)^{s-1} \leftrightarrow \Gamma(s) n^{-s}
\end{equation}

The integral of the $nth$ pole pair on the imaginary axis corresponds to the sum of the steps conjugate to step $n$. Riemann used the $\Gamma$-function identities to express this ratio in a form emphasizing the symmetry under $s \to 1-s$, as it is usually presented in the literature.

\begin{equation}
\zeta(s) = \Gamma (1-s) (2\pi)^{s-1} 2 sin(s\pi/2) \zeta (1-s) \qquad \sigma > 1.
\label{functeq}
\end{equation}

However, the ratio can also be expressed using Stirling's approximation $\Gamma(s) = (2\pi)^{1/2}(s-1)^{s-1/2} e^{1-s}\ldots$, giving $n_p^{1-2s}e^{i(\pi/4 + t + \ldots)}$ as derived above.  This form of $Q(s)$ emphasizes the symmetry of initial steps and conjugate regions, specifically displaying the dependence of $\lvert Q \rvert=n_p^{1-2\sigma}$ on $t$.

The functional equation of $\zeta$ is an immediate consequence of Equation \ref{sumconj}, except that this only describes ``half'' of the relation, the steps $n^{-s}$ and $n^{s-1}$ for $n \le n_p$.  Equation \ref{functeq} results from the sum of Equation \ref{sumconj} and its transform under $s \to 1-s$, each summed to $n_p$, and the relation $Q(s)Q(1-s)=1$.

\section{First-Order Determination of the Conjugate Points}

Figure \ref{xy9} strongly suggests that the point conjugate to the end of step $n$ is the center of the scroll having $\delta \theta = -(2n+1)\pi/2$, or $t = (2n+1)\pi$.  For this exact value of $t$, second order evaluation of $ \delta \theta_{N -1/2 \pm 1/2}$ gives the symmetric trio of steps $N-1,N,N+1$ with equally spaced angles $\pi \mp \pi/2N$.  When $t$ equals this odd multiple of $\pi$, the scroll center lies \textit{on} step $N$, giving a transverse displacement $\tau_N=0$. The midpoints of steps $N,N+1$ are separated by one half the step length $N^{-s}$ times their angular separation $\pi/2N$. Defining $\delta t = t- \pi[t/\pi]$, each increase of $\delta t$ by $\pi$ advances the scroll center to subsequent steps $N+1,N+2$, etc. Although the angle changes of the steps are large, equal to $\delta t log(N)$, the variation of the \textit{angle difference} $\delta \theta_N$ with $\delta t$ is small, yielding the linear relation $\tau_N = i\delta t/4N^{s+1}$.

For $\sigma \ne 0$, the variation of step length with $n$ leads to a longitudinal displacement $\lambda_N$ between the midpoint of a step and the geometric center of the scroll.  Defining that center as a fraction $\alpha$ of its step length, $\lambda_N$ is given by the condition that the centers of adjacent steps having angle changes of $\pi$ coincide.  Assuming that $\alpha$ varies weakly with $n$, $(1-\alpha)N^{-s}=\alpha(N+1)^{-s}$ gives $\alpha =1/2+\sigma/4N$.  The effect of this longitudinal displacement on $\tau_N$ is of higher order.

Taken together, and including the overall step angle $-tlog(n)$, the scroll center conjugate to step $n=t/\pi$ is given to first order by

\begin{equation}
\sum_1^N m^{-s} -\frac{1}{2N^{s}}+\frac{\sigma +i\delta t}{4N^{s+1}}+\ldots
\label{originprime}
\end{equation}

 The inclusion of higher order corrections for comparison with Riemann's analytic results is left to the reader.

\section{A Point Conjugate to Itself}

\begin{figure}[hp]
\centering
\includegraphics[angle=0,width=.5\textwidth]{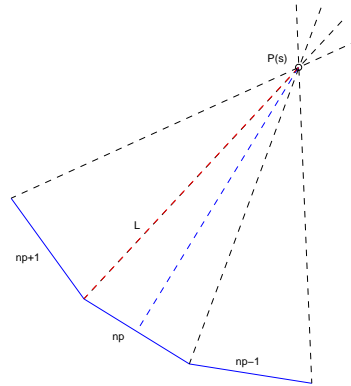}
\caption{Location of the Symmetry Center}
\label{centercnstr}
\end{figure}

The complete step plot has bilateral symmetry, but in general the point $n_p$ does not lie on the symmetry axis. Figure \ref{blob1pt9} clearly shows the ``pendant'' near step $n_p$ and its polygonal shape which approximates a circle for $n_p \gg 1$.  The center of this pendant is a point on the axis, denoted here by $P(s)$.  It is laterally displaced from $\sum_1^{n_p} n^{-s}$ by the distance $L$, as shown in Fig. \ref{centercnstr}.

Because $\delta^2 \theta = 2 \pi$ at step $n_p$, the steps $n_p$ and $n_p +1$ are tangent to this central feature.  Geometrical construction gives

\begin{equation}
P(s) = \sum_1^{n_p} n^{-s}+L
\label{center}
\end{equation}

\noindent where

\begin{equation}
L= -\frac{e^{-i(t log(n_p)+2 \pi p)}}{2 n_p^s \cos(2 \pi p)}.
\label{Loffset}
\end{equation}

\noindent  and where $p$ is the fractional part of $\sqrt{t/2\pi}$, given by $(n_p +p)^2 = t/2\pi$. Note that $L \ne 0$, as seen both by the construction and Equation \ref{Loffset}.  Figure \ref{centercnstr} also shows that $L$ is closely parallel to the axis of symmetry, as it is perpendicular to the average of angles $\theta_n$ and $\theta_{n+1}$.

\section{The Point Conjugate to the Origin}

The reciprocity relation between original and conjugate steps implies a final scroll near step $[t/\pi]$. Consider a pair of adjacent steps $n,n+1$ beyond the inflection point of the final conjugate region at $t/2\pi$, where $\delta \theta$ decreases uniformly to zero with no further discontinuities of $2\pi$ due to the mod function.  The perpendicular bisectors of these adjacent steps intersect at a distance $R = n^{-\sigma}/ \delta \theta \simeq n^{1-\sigma }/t$, defining a radius of curvature that decreases with increasing $n$ for $\sigma > 1$.  In this case, the final scroll collapses to the value of $\zeta$.  When $\sigma < 1$, $R$ increases without bound. \\

\textit{The series divergence is entirely due to the outward spiral of steps beyond the final scroll at $O'=\sum_1^{[t/\pi]} m^{-s}$ for $\sigma < 1$.} \\

Edwards notes that Riemann did not use the term ``analytic continuation'', but rather discussed ``a solution that remains valid'' when the series representation does not converge.  By that definition, the point $O'$ conjugate to the origin \textit{is} the analytic continuation of $\zeta(s)$ for $\sigma < 1$. This conclusion appears to differ from all methods reported previously, but is in agreement with them.  As one example, consider the \textit{Euler-Maclaurin summation method}, in which the explicit sum of terms beyond $n = [t/2\pi]$ is represented by an integration which extends the region of convergence to $\sigma > 0$.  This yields, to lowest order \cite[Sec.~6.4]{edwards}

\begin{equation}
\zeta (s) = \sum_1^N n^{-s} -\frac{N^{-s}}{2} + \frac{N^{1-s}}{s-1} \sum_{k=0} \frac{B_{2k} s^{2k}}{(2k)! N^{2k}}+\ldots
\end{equation}

\noindent where the \textit{Bernoulli numbers $B_n$} are given by the generating function

\begin{equation}
\frac{x}{e^x -1}= \sum_1^{\infty}\frac{B_n x^n}{n!}
\end{equation}

The sum over powers of $(s/N)$ converges slowly for $\lvert s/N \rvert >= \pi$, but it can be represented by the Bernoulli number generating function.   In this way, one obtains the first order correction of the final step;

\begin{equation}
\frac{N^{-s}}{2} \frac{s}{s-1} \frac{e^{s/2N}+e^{-s/2N}}{e^{s/2N}-e^{-s/2N}}
\end{equation}

The ratio $s/(s-1)$ approaches $1$ for large $s$, and manipulation of the definitions of trigonometric and hyperbolic functions lead to

\begin{equation}
\frac{N^{-s}}{2} \frac{\tanh(\sigma/2N)-i\cdot \cot(t/2N)}{1-i \cdot \tanh(\sigma/2N) \cot(t/2N)}
\end{equation}

 Expansion of the real term reproduces the longitudinal correction term $\sigma/4N^{s+1}$ for the final scroll. The term $\cot(t/2N)$ implements the mod function, and its expansion about the zero of the $cot$ gives $\delta t/4N^{s+1}$, resulting in agreement with the geometric first order evaluation of the point $O'$ conjugate to the origin.

\begin{equation}
\sum_1^{N=[t/pi]} n^{-s}-\frac{N^{-s}}{2}+\frac{\sigma + i\delta t}{4N^{s+1}}
\label{originprime}
\end{equation}

\section{First-Order Riemann-Siegel Correction}

Although Riemann did not include it in his paper, his handwritten notes (recovered by Siegel in a feat of exceptional mathematical ability) show that he used symmetry to facilitate calculations, including a hierarchy of equations giving successively higher order results. His first order result can be shown to correspond to $2L(s)$, the lowest order of Equation \ref{Loffset}.

In Riemann's development, $\sum_1^{n_p} n^{-s}$, its conjugate $Q(s) \sum_1^{n_p} n^{s-1}$ and a remainder term $R$ bridging the endpoints of the two sums are projected to the $\Theta$-axis on which $\zeta(s)$ resides, and termed $Z(s)$. ,

 \begin{equation}
 Z=(\sum_1^{n_p} n^{-s} + e^{2i\Theta}n_p^{1-2\sigma}\sum_1^{n_p}n^{s-1}+R(t))e^{-i\Theta}
 \label{cosarg}
\end{equation}

 Since $Z \in \Re$, locating its zeros is a matter of detecting sign changes. Riemann determined the first order ``mismatch'' between $\sum_1^{n_p} n^{-s}$ and its conjugate using the ``steepest descent'' or ``saddle-point'' method\cite[Sec.~7.2]{edwards}. The integrand of the mismatch is expanded about its peak. One segment of a closed contour crosses the imaginary axis at the $n_p th$ pole at the angle $\pi/4$, leading to an integral of gaussian form with its peak on the imaginary axis.  All other portions of this contour give negligible contributions. His first order expression, including his approximation for the remainder $R$ and assuming $\sigma = 1/2$, is

 \begin{multline}
 Z = \sum_1^{n_p} n^{-1/2} 2cos(\Theta(t)+t log(n))+ \\ (-1)^{n_p -1} \frac{\cos2\pi(p^2-p-\frac{1}{16})}{n_p^{1/2}\cos(2\pi p))}+\ldots
\label{bridge}
\end{multline}

 From symmetry, $L$ is related geometrically to Riemann's remainder $R$ by the expression $R=L(s)+Q(s)L(1-s)$. In the approximation that $L$ lies approximately along the symmetry axis this can be simplified to

 \begin{equation}
 R = 2\mathsf{L} \cos(\theta_L - \Theta)
 \end{equation}

 \noindent where

 \begin{multline}
\Theta = -t log(n_p) +\pi n_p^2 + \pi (\frac{\pi}{8}- 2 p^2) \\
\theta_L = -t log(n_p) - \frac{\pi}{2} -\frac{1}{2} \delta \theta \\
\end{multline}

 Using Equation \ref{Loffset} and noting that $\Theta$ involves $p$ through $t log(\sqrt(t/2\pi))=t log(n_p)+t log(1+p/n_p)$, discarding terms of $\mathcal{O}(1/n_p)$, and using $\pi n_p^2 (mod 2\pi)$ gives $0$ for $n_p$ even or $\pi$ for $n_p$ odd gives

\begin{equation}
\theta_L -\Theta = 2 \pi (p^2 - p -1/16 + (0  \quad or \quad \pi))
\end{equation}

\noindent Combined with the equation for $L$, this shows that Riemann's first order correction is given by the projection of $2L$ to the $\Theta$-axis. A study of Riemann's analysis in \cite[Sec.~7.4]{edwards} displays both his mathematical proficiency and the simplicity of the geometric approach.

The zeros in both the numerator and denominator of $R$ at $p=1/4,3/4$ can be understood geometrically. At these values of $p$, $\delta \theta_{n_p}$ is a multiple of $2\pi$ and steps $n_p$ and $n_{p+1}$ are each  exactly perpendicular to the symmetry axis.  The error in the pendant center location grows large at these values of $p$, but does so \textit{along} the axis.  It is therefore perpendicular to $\theta_L$.  In the projection to the $\Theta$-axis, the zeros of $cos(\Theta - \theta_L)$ cancel this effect.

\section{The ``Double Pendulum'' Construction}

\begin{figure}[pt]
\centering
\includegraphics[angle=0,width=.5\textwidth]{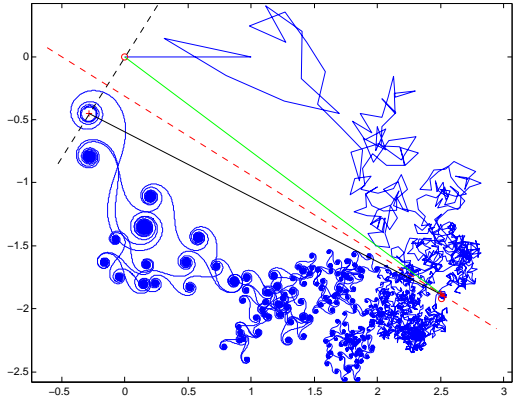}
\caption{Riemann-Siegel illustration for $s=1/2+2220000.15i$.}
\label{rslines}
\end{figure}

The Riemann-Siegel calculation is represented geometrically by the sum of two line segments, $OP=P(s)$ and $O'P = Q(s)P(1-s)$.  This ``double pendulum'', related to a classical geometrical figure called the \textit{lima\c con}, is illustrated in Figure \ref{rslines} for $s=1/2+2220000.15i$, chosen to reduce conjugate region overlaps and display the features clearly.  The green vector connects the origin to the pendant center $P(s)$.  This, plus its conjugate, in black, ends at $\zeta (s)$.  The bisector of these two lines is the symmetry axis, showing geometrically that $\zeta(s)$ lies on the line through the origin perpendicular to the symmetry axis.

\begin{figure}[pt]
\centering
\includegraphics[angle=0,width=.5\textwidth]{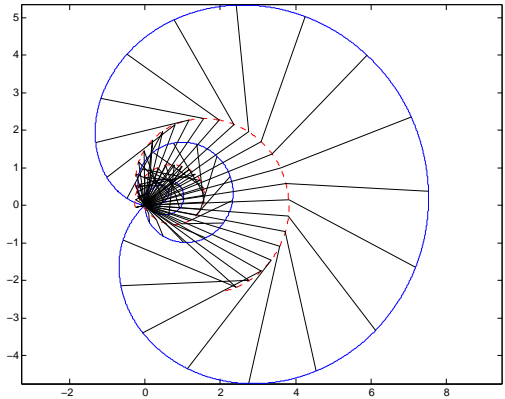}
\caption{Lima\c con construction for Gram points $1000 \to 1004$.}
\label{limacon1000}
\end{figure}

Figure \ref{limacon1000} illustrates this construction for a range of values of $t$ having $\sigma = 1/2$. Here, the dashed spiral shows $P(s)$ and the solid spiral gives $O' = \zeta(s)$ at successive values of $t$, with the curves connected at regularly spaced values of $t$ to indicate the line segments of the lima\c con construction.

\begin{figure}[pt]
\centering
\includegraphics[angle=0,width=.5\textwidth]{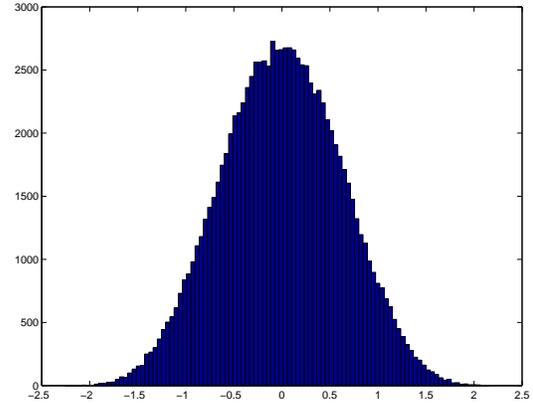}
\caption{Distribution of fractional positions of the zeros between adjacent Gram points, for the first 100,000 zeros.}
\label{fdist}
\end{figure}

 The number of zeros for $t < T$ has been a subject of much research.  Integration around the boundary of the critical strip shows that there are approximately as many zeros as the number of times that $\Theta(T)$ has passed through a multiple of $\pi$.  This leads to the definition of \textit{Gram points} $t_N$  by the relation $\Theta(t_N)=N\pi$ \cite[Sec.~6.5]{edwards}.  When $\Theta$ is an \textit{odd} multiple of $\pi /2$ the real part of $\zeta$ is zero, countering its ``bias'' to positive values because the first step always gives a real contribution of $+1$. In this case, nearby values of $t$ are more likely to lead to complex zeros.  Although the zeros are spaced irregularly because the angle of $P(s)$ does not change monotonically with $t$ as does $\Theta$, on average there \textit{does} appear to be one zero per Gram point and they \textit{are} found more frequently between successive Gram points. The tendency for zeros to ``group'' between Gram points is shown in Figure \ref{fdist}, a histogram of the displacement of the first 100,000 zeros \cite{odlyzko} from the midpoint of adjacent Gram points, scaled by half the distance between Gram points.

\begin{figure}[pt]
\centering
\includegraphics[angle=0,width=.45\textwidth]{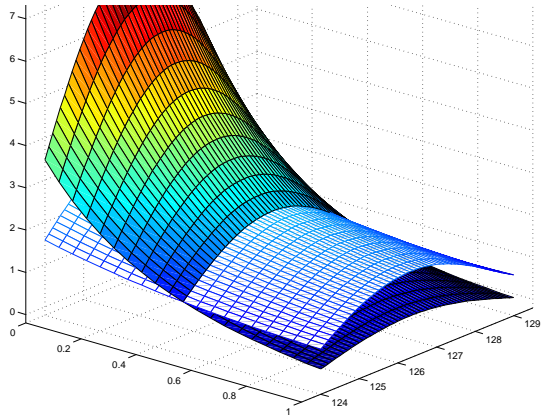}
\caption{Magnitudes of $P(s)$ (plain mesh) and $Q(s)P(1-s)$ (color) on the critical strip for $124<t<129$.}
\label{psurf120to129}
\end{figure}

When $\sigma = 1/2$, the condition of opposing angles is a one-parameter ``search'' as $t$ increases, just as the sign changes of $Z$ signal the crossing of a complex zero for $\zeta(t)$.  Thus, the symmetry shows why zeros occur frequently for this value of $\sigma$.  This has led to the \textit{Riemann hypothesis} \cite[Sec.~13.9]{apostol} that \textit{all} zeros are on this line in the complex plane.  Stated in this form, a proof of the hypothesis requires showing that zeros are impossible when $\sigma \ne 1/2$.  As illustrated in Figure \ref{psurf120to129}, $P(s)$ and $Q(s)P(1-s)$ for $120 < t < 129$ intersect transversely in this portion of the critical strip. This intersection is identical for $\sigma = 1/2$, where the transformation $s \to 1-s$ corresponds to complex conjugation and the magnitude of $Q(s)$ is unity.

\section{Zeros and the Critical Line}

 Riemann's \textit{On the Number of Prime Numbers Less Than a Given Quantity} was chosen for inclusion in \textit{God Created the Integers}\cite{hawking}, Hawking's list of landmark mathematical works. Riemann found that the first few zeros are on the critical line. However, this result was not obvious from the calculational algorithm; this peculiar generality was later stated as the general hypothesis bearing his name. Recognition of the symmetry sheds light on this, as can be illustrated by three specific algorithms for $\zeta$.  Since each is a solution of the same problem, presented from different viewpoints, they give the same numerical results within their respective error limits.

 The first algorithm evaluates $\zeta$ for a range of $t$ and arbitrary values of $\sigma$ using the Euler-Maclaurin approximation;

 \begin{equation}
 \zeta(s) \simeq \sum_1^{N=[t/\pi]} m^{-s} -\frac{1}{2N^{s}}+\frac{\sigma +i\delta t}{4N^{s+1}}
 \end{equation}

\begin{figure}[!]
\centering
\includegraphics[angle=0,width=.5\textwidth]{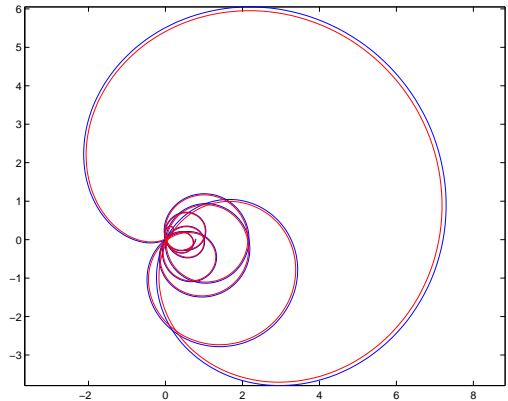}
\caption{Euler-Maclaurin, $2000<t<2010, blue:\sigma=.500,red:\sigma = .505$}
\label{aone}
\end{figure}

\begin{figure}[!]
\centering
\includegraphics[angle=0,width=.5\textwidth]{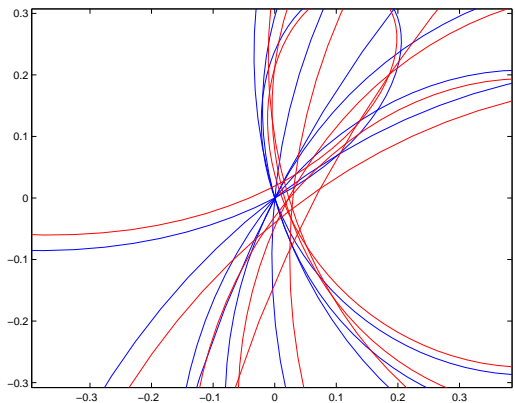}
\caption{Euler-Maclaurin enlarged, $2000<t<2010, blue:\sigma=.500,red:\sigma = .505$}
\label{aonelarge}
\end{figure}

 Figure \ref{aone} shows $\Im \zeta$ versus $\Re \zeta$ for $2000<t<2010$, with $\sigma=.500$ in blue and $\sigma = .505$ in red.  Enlargement of the region around the origin in Figure \ref{aonelarge} shows that every ``loop'' with $\sigma=1/2$ invariably yields a complex zero, whereas a slight variation of $\sigma$ gives ``loops'' which avoid the origin.

How can this extremely simple algorithm give a zero for every ``loop'', but only when $\sigma=1/2$?  When the author first encountered this result, having then no appreciation of the symmetry, it seemed wildly unusual, being nothing more than the sum of Dirichlet terms $n^{-s}$ with the first order Euler-Maclaurin correction.  What ``unseen hand'' could possibly force this particular result?  With symmetry, one realizes that although $P(s)$ and $O'$ are not explicitly represented in in this prescription, their effects are included correctly in the result. Note that the Euler-Maclaurin algorithm is insensitive to the upper limit $N>[t/2\pi]$; terminating the sum at $N+1$ and subtracting $(N+1)^{-s}/2$ gives the center of the next step, and the change of $\delta t$ by $\pi$ yields the same value for the spiral center as did the upper limit $N$.

Use of Riemann's formulation in Equation \ref{bridge} gives the same result far more efficiently, and in compact notation, summing only 17 terms instead of more than 600 for the Euler-Maclaurin algorithm with same range of values of $t$.

\begin{multline}
 Z = \sum_1^{n_p} n^{-1/2} 2cos(\Theta(t)+t log(n))+ \\ (-1)^{n_p -1} \frac{\cos2\pi(p^2-p-\frac{1}{16})}{n_p^{1/2}\cos(2\pi p))}
 \end{multline}

However, this form is limited to $\sigma = 1/2$ due to the compact expression of the projections of $\sum_1^{n_p}n^{-s}$ and $\sum_1^{n_p}n^{s-1}$ to the $\Theta$-axis using the cosine function. The remainder term is presented as a single entity bridging the two discrete sums, masking the identity of the symmetry center.

\begin{figure}[!]
\centering
\includegraphics[angle=0,width=.5\textwidth]{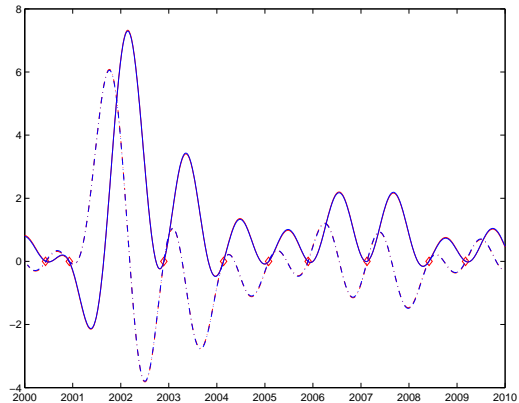}
\caption{$\zeta$ Real (-) and imaginary (- .) components vs. $t$ for Riemann's form (blue) and the symmetric form (red).  Red diamonds are the actual zeros.}
\label{symmrs}
\end{figure}

Implementation of the symmetric form

 \begin{equation}
 \zeta(s) = P(s)+Q(s)P(1-s)
 \end{equation}

  \noindent using Equations \ref{center} and \ref{Loffset} demonstrates the symmetry explicitly and is valid for arbitrary $\sigma$. It gives the results and numerical efficiency of Riemann's form, as shown in Figure \ref{symmrs}. Seen from this perspective, however, there is no peculiarity in the preference for zeros when $\sigma=1/2$.

The functional equation demands that zeros off of the critical line must appear in pairs at $s = 1/2 \pm \epsilon + it$, and exceptional conditions are required for this to occur.  The symmetric Riemann-Siegel equation shows that either:

\begin{quote}\textit{The magnitudes of $P(s)$ and $Q(s)P(1-s)$ must be equal for zeros at $\sigma = 1/2 \pm \epsilon$  $\textbf{and}$ the condition of opposing angles must be met,}\end{quote}

 or:

\begin{quote}\textit{Since $Q(s) \ne 0$ on the critical strip, $P(1/2+\epsilon+it)$ and $P(1/2-\epsilon+it)$ must be simultaneously zero at these arguments.}\end{quote}

Symmetry, with the example calculation shown in Figure \ref{psurf120to129}, exposes the difficulty of attaining these conditions; it may, however, be the case that they cannot be excluded ``in the fullness of $t$''.  If so, this would explain why the Riemann Hypothesis has resisted any proof for over 150 years\cite{wolf}.

\section{Acknowledgements}
Interest in this subject was kindled many years ago by a statistical physics professor who chose \emph{not} to say "This integral is $\pi^4 /15$, you can look it up", but gave the class a lecture on the zeta function, including ``the sport that mathematicians call analytical continuation''.  Michael R. Stamm of the Air Force Research Laboratory has provided the most consistent help and encouragement.  Mark Coffey of the Colorado School of Mines has shared many preprints and offered references for $\zeta$ function applications.  Steve Lamoreaux of the Physics Department at Yale University has suggested that Riemann's analysis may have been guided by geometric reasoning. The assistance of Billy Buttler of the Physics Division at the Los Alamos National Laboratory not only for mathematical advice, but for calculations and figures with Matlab has been absolutely indispensable.

\end{document}